\documentclass[a4paper,12pt]{amsart}
\usepackage{graphicx}
\usepackage{color}

\numberwithin{equation}{section}
\numberwithin{figure}{section}

\newtheorem{thm}{Theorem}
  
\newtheorem{thrm}{Theorem}[section]
\newtheorem{prop}[thrm]{Proposition}
\newtheorem{lemma}[thrm]{Lemma}

\theoremstyle{definition}
 \newtheorem*{dfn}{Definition}
 
 \newtheorem{expl}{Example}%

\theoremstyle{remark}
 \newtheorem*{rmk}{Remark}

\renewcommand{\phi}{\varphi}
\renewcommand{\epsilon}{\varepsilon}
\newcommand{\rgt}[1]{\right#1}
\newcommand{\lft}[1]{\left#1}
\newcommand{\rd}{\partial}
\newcommand{\mbdr}{\partial_-}
\newcommand{\rhdl}[2][k]{D^{#1}\times D^{#2}\times S^1}
\newcommand{\R}{\mathbb{R}}
\newcommand{\Z}{\mathbb{Z}}
\newcommand{\inter}{\operatorname{int}}

\newcommand{\im}{\operatorname{Im}}
\newcommand{\releu}[2][\xi]{\langle e(#1),[#2]\rangle}
\newcommand{\dsc}[1]{D(#1)}
\newcommand{\stm}{S^1\times}

\begin{document} 
%
%
\title[Round surgery and contact structures on $3$-manifolds]
{Round surgery and \\ contact structures on $3$-manifolds}
\author{Jiro ADACHI}
\thanks{This work was supported 
  by JSPS KAKENHI Grant Numbers~21540058, 25400077, 
  and MEXT KAKENHI Grant Number~17075010
. }
\address{Department of Mathematics, 
  Hokkaido University, 
  Sapporo 060--0810, Japan.}
\email{j-adachi@math.sci.hokudai.ac.jp}

\begin{abstract}
  Contact round surgery of contact $3$-manifolds is introduced in this paper. 
 By using this method, 
an alternative proof of the existence of a contact structure 
on any closed orientable $3$-manifold is given. 
 It is also proved 
that any contact structure on any closed orientable $3$-manifold 
is constructed from the standard contact structure 
on the $3$-dimensional sphere by contact round surgeries. 
 For the proof, important operations in contact topology, 
the Lutz twist and the Giroux torsion, are described 
in terms of contact round surgeries. 
\end{abstract}

\maketitle

\section{Introduction}\label{sec:intro}
  Surgery is a basic and important method in studying topology of manifolds. 
 Not only for manifolds themselves 
but also for geometric structures on manifolds, surgery is important. 
 In this paper, a new notion of round surgery with contact structures 
is introduced. 
 It is proved that all closed connected contact $3$-manifolds are constructed 
by this method. 
 It would give a new perspective on the study of contact structures. 

  Round surgery of a manifold was originally introduced by Asimov~\cite{asimov}.
 He introduced the notion of round handle 
in order to study non-singular Morse Smale flows 
(see~\cite{asimov}, \cite{morgan}). 
 Further, the theory had been applied to some geometric studies 
(see~\cite{miyoshi}, \cite{etgh}, \cite{vogel}, \cite{baykur}). 
 A \emph{round handle}\/ is roughly defined 
as a product of an ordinary handle with a circle. 
 In other words, a round handle of index~$k$ of dimension $n$ is 
$R_k=D^k\times D^{n-k-1}\times S^1$ 
attached to the boundary $\rd N$ of an $n$-dimensional manifold $N$ 
by an embedding $\phi\colon\rd D^k\times D^{n-k-1}\times S^1\to\rd N$. 
 Round surgeries of an $m$-dimensional manifold $M$ is defined 
through cobordisms
by attaching round handles of dimension $m+1$ to $M\times[0,1]$ 
(see Section~\ref{sec:rd-handle} for precise definition). 
 In this paper, we are mainly dealing with connected $3$-dimensional manifolds. 
 Therefore, there exist two kinds of round surgeries, index~$1$ and index~$2$. 
 Asimov proved in \cite{asimov} that any closed orientable $3$-manifold 
is obtained from a $3$-dimensional sphere by a sequence of round surgeries 
of index~$1$ and index~$2$. 

  One of the purposes of this paper is to apply this surgery theory 
to contact $3$-manifolds. 
 A \emph{contact structure}\/ on a $3$-dimensional manifold 
is a completely non-integrable plane field 
(see Subsection~\ref{sec:cttop} for precise definition). 
 It is known that on any closed orientable $3$-manifold, 
there exists a contact structure. 
 It is proved firstly by Martinet (\cite{martinet}, \cite{lutz}), 
and there are some alternative proofs (see \cite{geitext}). 
%
%
\begin{thm}\label{thma} 
  On any closed orientable $3$-manifold, there exists a contact structure. 
\end{thm} 
\noindent
 One of the results in this paper is 
an alternative proof of this existence theorem. 
 By defining contact round surgery, 
that is, round surgery with contact structure, 
we apply Asimov's result above. 
 By contact round surgeries, a contact structure is constructed 
on any closed orientable $3$-manifold (see Subsection~\ref{sec:prThma}). 

  Not only the existence of a contact structure, 
we show that all contact structures on all closed orientable $3$-manifolds 
are constructed by contact round surgeries 
from the standard contact $3$-dimensional sphere. 
%
%
\begin{thm}\label{thmb}
  Any contact structure on any closed orientable $3$-manifold 
is constructed from the standard contact $3$-dimensional sphere 
by a sequence of contact round surgeries of index~$1$ and index~$2$. 
\end{thm} 
\noindent
 In the proof of this theorem, all closed connected contact $3$-manifolds 
are reduced to contact structures on $S^3$ by contact round surgeries 
(Subsection~\ref{sec:prThmb}). 
 Then the classification of contact structures on $S^3$ 
due to Eliashberg~(\cite{eliashOT}, \cite{eliash20years}) is applied. 
 It is known that contact structures on $S^3$ are described by the Lutz twists, 
a modification of contact structures introduced by Lutz~\cite{lutz}. 
 In order to apply the classification result, 
it is proved that the Lutz twist is represented 
by a certain sequence of contact round surgeries
(see Theorem~\ref{thm:Lz}).

  It should be remarked that the realization of the generalized version 
of the Lutz twist along a certain torus is given in this paper. 
 It is a distinctive point of this method 
compared with the Weinstein surgery (see \cite{weinsteinHdl}, \cite{dingei04}). 
 Because of this, as well as the Lutz twist, 
it is proved that the Giroux torsion is closely related 
to contact round surgeries (see Subsection~\ref{sec:dfnlztw}). 
 The Giroux torsion is a notion introduced by Giroux~\cite{girouxTor}, 
which is a source of symplectically non-fillable but tight contact structures. 

  Some kinds of contact surgeries had appeared 
in the study of contact structure. 
 We should first mention contact Dehn surgery along transverse links
due to Lutz~\cite{lutz} and Martinet~\cite{martinet}. 
 By using this surgery method, they proved the existence of a contact structure 
on any closed orientable $3$-manifold (Theorem~\ref{thma}). 

  Symplectic or Stein handlebody surgery was studied 
by Weinstein~\cite{weinsteinHdl} and Eliashberg~\cite{eliashCharStein}. 
 It is an important tool to study the fillability of contact manifolds 
by symplectic or Stein manifolds. 

  Contact Dehn surgery along Legendrian link is now actively studied 
(see for example \cite{osbook}). 
 In terms of this surgery, the symplectic handlebody surgery is written as 
$(-1)$-surgery. 
 Note that the surgery coefficients in this case are measured with respect to 
contact framings of Legendrian knots. 
 Ding and Geiges~\cite{dingei04} proved that any rational contact Dehn surgery 
along a Legendrian knot is represented by contact Dehn surgeries 
of surgery coefficients $\pm1$. 
 Further, they proved that any contact structure 
on any closed orientable $3$-manifold 
is constructed by contact $(\pm1)$-surgeries. 
 Etnyre~\cite{etsurg} illuminated the study of $(+1)$-surgeries. 
 After  Stipsicz~\cite{stipsicz} related these surgeries 
to open book decomposition, 
these surgeries are studied with the Heegaard-Floer theory 
and the Ozsv\'ath-Szab\'o contact invariants (see \cite{ozsz}). 

  This article is organized as follows. 
 In the next section, we review some basic notions on contact topology 
which are needed in the discussions later. 
 First, we review some basic definitions and properties, 
and then some useful results in the convex surface theory. 
 In Section~\ref{sec:rd-handle}, we review the round handle theory. 
 Contact round surgeries are defined in Section~\ref{sec:c-rd-surg}. 
 In Section~\ref{sec:LtwGtor}, a relationship between contact round surgeries 
and the Lutz twist is discussed. 
 We also mention a relationship with the Giroux torsion. 
 Finally, Theorem~\ref{thma} and Theorem~\ref{thmb} are proved 
in Section~\ref{sec:concl}. 

  Several years has passed since the former version was written. 
 Meanwhile, the results in this paper has been considerably improved. 
 Another construction of the Lutz twist is given in~\cite{art20}. 
 Furthermore, symplectic round handles with the Liouville vector fields 
corresponding to contact round surgeries 
are also constructed by the author in any even dimension~\cite{art19}. 
 Researches on contact structures on higher-dimensional manifolds 
have advanced (see~\cite{mnw}, \cite{bem}). 
 Applications of contact round surgery to higher-dimensions are expected. 

\smallskip

  The author obtained the idea of this work when he visited
Technion -- Israel Institute of Technology. 
 He would like to thank Professor Michail Zhitomirski\u{\i} 
for the warm hospitality and especially for giving him a chance. 
 The author is also grateful to Professor Yakov Eliashberg, 
Professor John Etnyre, and Professor Ko Honda 
for some valuable discussions. 

\section{Contact geometry}\label{sec:prelim}
  In this section, we introduce basic notions and techniques 
in contact geometry 
which are needed in the discussion later. 
 First, we introduce some basic notions. 
 And then, we review convex surface theory 
which is a useful technique in contact topology. 

\subsection{Contact structures on $3$-manifolds} 
\label{sec:cttop}

  Let us begin with basic definitions. 
 A \emph{contact structure}\/ on a $3$-dimensional manifold $M$ 
is defined to be a tangent plane field $\xi$ on $M$ 
which is completely non-integrable. 
 In other words, contact structure $\xi$ is represented locally as the kernel 
$\xi=\ker\alpha$ of a $1$-form $\alpha$ 
which satisfies $\alpha\wedge d\alpha\ne 0$. 
 Note that the sign of the $3$-form $\alpha\wedge d\alpha$ does not depend 
on the choice of the $1$-form $\alpha$ but on $\xi$ itself. 
 Therefore, a $3$-manifold with a contact structure is orientable. 
 In this paper, we assume that $3$-manifolds are oriented, 
and that contact structures are \emph{positive}, that is, 
each corresponding $3$-form is a positive volume form 
on each oriented manifolds. 

  Some of the most basic properties of contact structures 
are the local triviality and the global stability. 
 Any contact structure $\xi$ on a $(2n+1)$-dimensional manifold 
is locally equivalent to the standard contact structure 
$\xi_0:=\lft\{dz+\sum_{i=1}^n(-y_idx_i+x_idy_i)=0\rgt\}$ on $\R^{2n+1}$, that is, 
there exists a local diffeomorphism which maps $\xi$ to $\xi_0$ 
(the Darboux theorem). 
 This implies that there is no local invariant for contact structures. 
 On the other hand, a deformation of a contact structure on a closed manifold 
can be traced by a one-parameter family of global diffeomorphisms 
(the Gray theorem). 
 This implies that contact structures are flexible. 
 Here, we should mention the local triviality of contact structures 
along curves transverse to the structures due to Lutz and Martinet. 
 Let $\xi_0$ be the standard contact structure on $S^1\times\R^{2n}$ 
defined by the $1$-form $d\phi+\sum_{i=1}^n(-y_idx_i+x_idy_i)$, 
where $\phi$ is a coordinate of $S^1$ and $(x_1,y_1,\dots,x_n,y_n)$ are 
those of $\R^{2n}$. 

%
%
\begin{thrm}[Lutz~\cite{lutz}, Martinet~\cite{martinet}]\label{thm:s1darboux}
  Let $\gamma$ be an embedded circle in a $(2n+1)$-dimensional manifold $M$, 
and $\xi$ a contact structure defined on a tubular neighborhood of $\gamma$. 
 Assume that $\xi$ is transverse to $\gamma$ at any point of $\gamma$. 
 Then there exists a local diffeomorphism 
from a tubular neighborhood of $\gamma$ 
to a tubular neighborhood of $S^1\times\{0\}\subset S^1\times\R^2$ 
which maps $\gamma$ to $S^1\times\{0\}$ and maps $\xi$ to $\xi_0$. 
\end{thrm}

\noindent
 In contact topology of $3$-dimensional manifolds, 
it is important to know contact structures on neighborhoods 
of embedded surfaces.
 We will discuss it later in the next subsection. 

  On curves transverse to contact structures, 
the following property is well-known (see for example~\cite{geitext}). 
%
%
\begin{prop} 
\label{prop:trsvApp}
  For any curve $L$ in a contact $3$-manifold, 
there exist positively and negatively transverse curves $L_\pm$ 
which are $C^0$-close to $L$. 
\end{prop}

  It is well-known that contact structures on $3$-dimensional manifolds 
are divided into two contradictory classes, tight and overtwisted. 
 A contact structure $\xi$ on a $3$-dimensional manifold $M$ 
is said to be \emph{overtwisted\/} 
if there exists an embedded disk $D\subset M$ which is tangent to $\xi$ 
along its boundary: $T_xD=\xi_x$ at any point $x\in \rd D$. 
 A contact structure $\xi$ is said to be \emph{tight \/} 
if it is not overtwisted. 

  Classification of contact structures is an important problem. 
 The classification, up to isotopy, of overtwisted contact structures 
is obtained by Eliashberg~\cite{eliashOT}. 
 It is proved that overtwisted contact structures 
which are homotopic as plane fields on a $3$-dimensional manifold are isotopic.
 In other words, the classification is reduced to some algebraic arguments 
in this case. 
 Concerning tight contact structures, 
although classifications on some manifolds had been obtained, 
it is still open on many manifolds. 
 In this paper, we need the classification of contact structures 
on the $3$-dimensional sphere $S^3$. 
 In this case, it is proved by Eliashberg~\cite{eliash20years} 
that there exists a unique tight contact structure called 
the standard contact structure on $S^3$. 
 Combining these results, 
we have the complete classification of contact structures on $S^3$. 
 With some trivialization, the homotopy classes of tangent plane fields 
are identified with the class of $\pi_3(S^2)\cong\Z$, 
by taking the Gauss mapping. 
 The classification is as follows. 

%
%
\begin{thrm}[Eliashberg, \cite{eliash20years}]\label{thm:classifS3}
  One homotopy class in $\pi_3(S^2)\cong\Z$ 
contains two non-equivalent contact structures on $S^3$, tight and overtwisted. 
 Any other class in $\pi_3(S^2)\cong\Z$ 
contains a unique overtwisted contact structure on $S^3$. 
\end{thrm}

\subsection{Convex surface theory}\label{sec:convsurf}
  Convex surface theory is one of key tools 
which caused a breakthrough on $3$-dimensional contact topology. 
 We review some basic properties of convex surfaces in contact $3$-manifolds 
which are needed in the definition of contact round surgery 
and the discussion below. 

  Before defining convexity of a surface in a contact $3$-manifold, 
we observe $1$-dimensional singular foliation on the surface 
traced by the contact plane field. 
 Let $F$ be an embedded surface in a contact $3$-manifold $(M,\xi)$. 
 The contact structure $\xi$ traces a singular $1$-dimensional foliation on $F$.
 In other words, $\xi_x\cap T_xF$ is a line in $T_xM$ 
if $\xi_x$ and $T_xF$ do not coincide at $x\in M$. 
 Then we obtain a line field ${\mathcal F}$ with singularities 
where $\xi_x=T_xF$. 
 By integrating the line field, 
we obtain a $1$-dimensional foliation with singularities on $F$. 
 Such a foliation is called the \emph{characteristic foliation} on $F$ 
with respect to $\xi$. 
 Let $F_\xi$ denote it. 
 It is known that a characteristic foliation on a surface 
determines a germ of contact structures along the surface. 

%
%
\begin{prop}[Giroux, \cite{girouxConv}]\label{prop:chfol}
  Let $F$ be an embedded surface in each contact $3$-manifold $(M_i,\xi_i)$, 
$i=1,2$. 
 Assume that they have the same characteristic foliation\textup{:} 
$F_{\xi_1}=F_{\xi_2}$. 
 Then there exists a local diffeomorphism 
between neighborhoods of $F\subset M_1$ and $F\subset M_2$
mapping the contact structure $\xi_1$ to $\xi_2$ and the surface $F$ to $F$. 
\end{prop}

\noindent
 This implies that two contact manifolds with boundaries 
can be glued together if they have the same characteristic foliation 
with orientation on their boundaries. 

  Now we define convexity. 
 It is defined by using the following vector field. 
 A vector field $X$ on a contact $3$-manifold $(M,\xi)$ 
is said to be \emph{contact\/} 
if its flow $\phi_t$ preserves the contact structure $\xi$: 
$\lft({\phi_t}\rgt)_\ast\xi=\xi$. 
 Let $F$ be a surface embedded in a contact $3$-manifold $(M,\xi)$. 
 The surface $F$ is said to be \emph{convex\/} 
if there exists a contact vector field on a neighborhood of $F\subset M$ 
which is transverse to $F$. 
 It is known that any surface in a contact $3$-manifold 
can be approximated to a convex surface. 

%
%
\begin{thrm}[Giroux, \cite{girouxConv}]\label{thm:convGen}
  For any embedded surface $F$ in a contact $3$-manifold $(M,\xi)$, 
there exists a convex surface $\tilde F\subset(M,\xi)$ 
which is $C^\infty$-close to $F$. 
\end{thrm}

  The most important reason why we adopt convex surfaces is a flexibility. 
 In order to describe the property, we need the following notions. 
 Let $\Sigma$ be a convex surface in a contact $3$-manifold $(M,\xi)$, 
and $X$ a contact vector field defined near $\Sigma$ 
which is transverse to $\Sigma$. 
 The \emph{dividing set}\/ $\Gamma_\Sigma$ of $\Sigma$ is defined as 
\begin{equation*}
  \Gamma_\Sigma:=\{x\in\Sigma\mid\xi_x\ni X_x\}. 
\end{equation*}
 In other words, it is a set of points where the contact plane gets ``vertical''
to $\Sigma$. 
 Note that $\Gamma_\Sigma$ does not depend on the choice of $X$ up to isotopy. 
 It is known that a dividing set is a curves 
transverse to leaves of a characteristic foliation (\cite{girouxConv}). 
 On the other hand, 
a $1$-dimensional foliation $\mathcal{F}$ on $\Sigma$ with singularity
is said to be \emph{divided\/} by $\Gamma_\Sigma$ if
\begin{itemize}
\item $\Sigma\setminus\Gamma_\Sigma=U_+\sqcup U_-$, 
\item $\Gamma_\Sigma$ is transverse to any leaf of $\mathcal{F}$, 
\item there exists a vector field $v$ which is tangent to $\mathcal{F}$ 
  and a volume form $\omega$ on $\Sigma$ which satisfy: 
  \begin{enumerate}
  \item $\operatorname{div}_\omega X>0$ on $U_+$, 
    $\operatorname{div}_\omega X<0$ on $U_-$, 
  \item the vector field $v$ looks outward of $U_+$ at $\Gamma_\Sigma$. 
  \end{enumerate}
\end{itemize}
 The following theorem makes dividing sets more flexible 
than characteristic foliations. 

%
%
\begin{thrm}[Giroux, \cite{girouxConv}]\label{thm:convFlex}
  Assume that $\Sigma$ is a convex surface in a contact $3$-manifold $(M,\xi)$ 
with a contact vector field $X$ transverse to $\Sigma$, 
and that $\Gamma_\Sigma\subset\Sigma$ is a dividing set. 
 Let $\mathcal{F}$ be a $1$-dimensional foliation with singularity on $\Sigma$ 
divided by $\Gamma_\Sigma$. 
 Then there exists a family $\phi_t\colon\Sigma\to M$, $t\in[0,1]$, 
of embeddings which satisfies\textup{:} 
\begin{itemize}
\item $\phi_0=\textup{id}_\Sigma$, $\phi_t|_{\Gamma_\Sigma}
  =\textup{id}_{\Gamma_\Sigma}$, 
  for any $t\in[0,1]$,
\item $\phi_t(\Sigma)$ is transverse to $X$ for any $t\in[0,1]$,
\item $(\phi_1(\Sigma))_\xi=(\phi_1)_\ast\mathcal{F}$. 
\end{itemize}
\end{thrm}

\noindent
 This Theorem~\ref{thm:convFlex} with Theorem~\ref{thm:convGen} 
implies that it is sufficient to check dividing sets on the boundaries 
and their orientations 
in order to glue two contact manifolds. 

\begin{figure}[htb]
  \centering
  \def\svgwidth{10.8cm}
  {\small 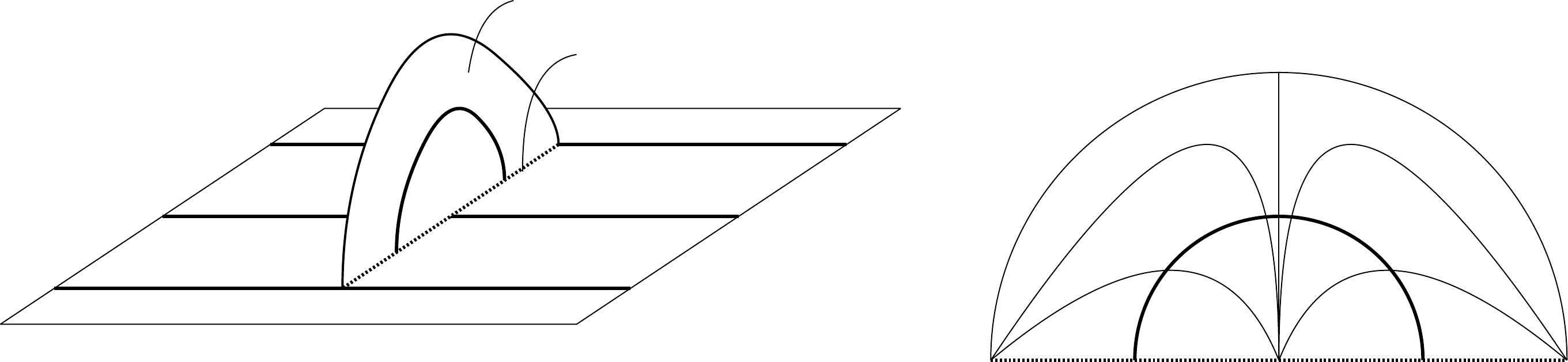}
  \caption{bypass}
  \label{fig:bypass}
\end{figure}

  There is a remarkable method to configure dividing curves 
developed by Honda~\cite{hondaI}. 
 A notion called bypass is introduced as follows. 
 Let $\Sigma$ be a convex surface in a contact $3$-manifold, 
$\Gamma_\Sigma\subset\Sigma$ a dividing set, 
and $\alpha\subset\Sigma$ a Legendrian arc 
that intersects $\Gamma_\Sigma$ transversely in three points $p_1,\ p_2,\ p_3$, 
where $p_1,\ p_3$ are end points of $\alpha$. 
 A \emph{Legendrian}\/ curve is a curve in a contact $3$-manifold 
which is everywhere tangent to the contact plane field. 
 A convex half-disk $D$ with Legendrian boundary 
is called a \emph{bypass}\/ for $\Sigma$ along $\alpha$ 
if it satisfies the following conditions (see Figure~\ref{fig:bypass}): 
\begin{itemize}
\item $D$ intersects $\Sigma$ along $\alpha$ transversely 
  on its boundary $\rd D$, 
\item $\operatorname{tb}(\rd D)=-1$. 
\end{itemize}
 The Thurston-Bennequin invariant $\operatorname{tb}(\gamma)$ 
of a closed Legendrian curve $\gamma$ 
is a twisting number of the contact plane field along the curve 
with respect to a certain framing. 
 In this case the framing depends on the disk $D$. 
%
%
 The standard characteristic foliation on a bypass half-disk appears 
as in Figure~\ref{fig:bypass}, where the thick curve is the dividing curve. 
 By attaching a bypass, dividing set is configured as follows. 
%
%
\begin{lemma}[Honda, \cite{hondaI}]\label{lem:bpatt}
  Let $\Sigma$ be a convex surface in a contact $3$-manifold. 
 Assume that there exists a bypass $D$ for $\Sigma$ 
along a Legendrian arc $\alpha$. 
  Then there exists a neighborhood of $\Sigma\cup D$ 
diffeomorphic to $\Sigma\times[0,1]$ 
which satisfies the following conditions\textup{:} 
\begin{itemize}
\item $\Sigma=\Sigma\times\{\epsilon\}$ for some $\epsilon\in(0,1)$, 
\item $\Sigma\times[0,\epsilon]$ is invariant, 
  that is, defined by a contact vector field transverse to $\Sigma$, 
\item $\Sigma\times\{1\}$ is convex, 
\item the dividing curve $\Gamma_{\Sigma\times\{1\}}$ is obtained 
    from $\Gamma_\Sigma$ by the operation 
    in \textup{Figure}~\textup{\ref{fig:bypatt}} 
    in a neighborhood of $\alpha$. 
\end{itemize}
%
%
\begin{figure}[htb]
  \centering
  \def\svgwidth{10.8cm}
  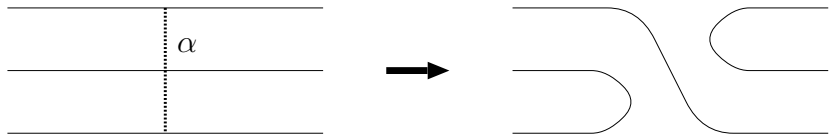
  \caption{attaching a bypass}
  \label{fig:bypatt}
\end{figure}
\end{lemma} 

  Finding an embedded bypass, we obtain a perturbation of a surface 
that causes the change of dividing curves by an attachment of the bypass 
like Figure~\ref{fig:bypatt}. 
 In general it is not easy to find a suitable embedded bypass. 
 When a contact structure is overtwisted, we can find a bypass we need. 

%
%
\begin{prop}[Huang, \cite{huang}]\label{prop:huang}
  Let $\Sigma$ be a convex surface in a contact $3$-manifold $(M,\xi)$. 
 Assume that $\xi$ is overtwisted on $M\setminus\Sigma$. 
 Then, for any Legendrian arc $\alpha$ on $\Sigma$ 
as the definition of the bypass, 
there exist bypasses along $\alpha$ in $M\setminus\Sigma$ 
attached from the both sides of $\Sigma$. 
\end{prop} 

  A criterion for overtwistedness from dividing curves is introduced by Giroux. 
%
%
\begin{prop}[Giroux,\cite{girouxCrit}]\label{prop:girouxcrit}
  Let $\Sigma$ be a closed convex surface in a contact $3$-manifold $(M,\xi)$. 
 Unless $\Sigma=S^2$ and the dividing set $\Gamma_\Sigma$ on $\Sigma$ 
is connected, 
the invariant tubular neighborhood of $\Sigma$ is overtwisted 
if $\Gamma_\Sigma$ contains a circle that is contractible in $\Sigma$. 
\end{prop} 

\section{Round handle theory}\label{sec:rd-handle}
  Round handle theory is introduced by Asimov~\cite{asimov} 
to study the Morse-Smale flow. 
 In this paper, we need an application to $3$-manifold theory by himself. 

  Round handle and round handle decomposition are defined as follows. 
 Let $M$ be a manifold of dimension $n$ with boundary $\rd M\ne\emptyset$. 
%
%
\begin{dfn}
  A \emph{round handle}\/ of dimension $n$ and index~$k$ attached to $M$ 
is defied as a pair
\begin{equation*}
  R_k=\lft(\rhdl{n-k-1}, f\rgt)
\end{equation*} 
of a product of an $(n-1)$-dimensional disk $D^k\times D^{n-k-1}$ with a circle
and an attaching embedding 
$f\colon\mbdr\lft(\rhdl{n-k-1}\rgt)\to\rd M$, 
where $\mbdr\lft(\rhdl{n-k-1}\rgt):=\rd D^k\times D^{n-k-1}\times S^1$ 
is the attaching region. 
 Let $M\cup_fR_k$ or $M+R_k$ denote the manifold obtained 
from $M$ and $\rhdl{n-k-1}$ by the attaching mapping $f$. 
\end{dfn} 
\noindent
 Sometimes, $R_k$ also denotes $\rhdl{n-k-1}$ itself 
or the corresponding subset in $M\cup_fR_k$. 
 A manifold $M$ is said to have a \emph{round handle decomposition}\/ 
if it is obtained from $N\times I$ by attaching round handles: 
\begin{equation*}
  M=(N\times I)+R_0^1+\dots+R_0^{i_0}+\dots+R_{n-1}^{i_{n-1}}, 
\end{equation*}
where $N$ is an $(n-1)$-dimensional manifold without boundary 
and $R_k^i$ are round handles of index~$k$. 

  Round handles are used to study flow manifolds. 
 A flow manifold is defined as follows. 
 Let $(M,\mbdr M)$ be a pair of a manifold $M$ with a specific union $\mbdr M$ 
of connected components of the boundary $\rd M$. 
 The pair $(M,\mbdr M)$ is called a \emph{flow manifold}\/ 
if there exists a non-singular vector field on $M$ 
which looks inward on $\mbdr M$ and outward on $\rd_+M:=\rd M\setminus\mbdr M$. 
 The following property of flow manifolds is proved by Asimov. 

%
%
\begin{thrm}[Asimov, \cite{asimov}] \label{thm:rh-decomp}
  Let $(M,\mbdr M)$ be a compact flow manifold 
whose dimension is greater than $3$. 
 Then, $M$ has a round handle decomposition. 
\end{thrm}

  By defining round surgery, 
the result above is applied to the study of $3$-dimensional manifolds. 
 The surgery is defined by using round handles in stead of ordinary handles. 
 In other words, a round surgery corresponds to attaching a round handle 
to a cobordism. 
 Let $M$ be a manifold of dimension $n$. 
 A \emph{round surgery}\/ of index~$k$ is defined as the operation 
removing an embedded $\inter\lft(\rd\rhdl{n-k}\rgt)$ from $M$ 
and regluing $D^k\times\rd D^{n-k}\times S^1$ by the identity mapping of 
$\rd D^k\times\rd D^{n-k}\times S^1$. 
 Applying Theorem~\ref{thm:rh-decomp} 
to a cobordism between two $3$-dimensional manifolds, 
Asimov proved the following theorem: 

%
%
\begin{thrm}[Asimov, \cite{asimov}] \label{thm:r-surg}
  Let $M$ be a connected closed orientable manifold of dimension $3$. 
 Then $M$ can be obtained from a $3$-dimensional sphere $S^3$ 
by a finite sequence of round surgeries of index~$1$ and index~$2$. 
\end{thrm}

  In the $3$-dimensional case, round surgeries of index~$1$ and index~$2$ are 
explicitly described as follows. 
 A round surgery of index~$1$ is the operation 
removing two open solid tori 
$\inter\lft(\rd D^1\times D^2\times S^1\rgt)
=\{\text{two points}\}\times(\inter D^2\times S^1)$ 
from a $3$-manifold $M$ 
and regluing a thickened torus $D^1\times\rd D^2\times S^1=I\times T^2$
by the identity mapping of a pair of tori 
$\rd D^1\times\rd D^2\times S^1=\{\text{two points}\}\times T^2$. 
 A round surgery of index~$2$ is the operation 
removing an open thickened torus 
$\inter\lft(D^1\times\rd D^2\times S^1\rgt)=\inter\lft(I\times T^2\rgt)$
from a $3$-manifold $M$ 
and regluing two solid tori 
$\rd D^1\times D^2\times S^1=\{\text{two points}\}\times(D^2\times S^1)$ 
by the identity mapping of a pair of tori 
$\rd D^1\times\rd D^2\times S^1=\{\text{two points}\}\times T^2$. 

\section{Contact round surgery} \label{sec:c-rd-surg}
  Contact round surgeries of a contact $3$-manifold 
of index~$1$ and index~$2$ are defined in this section. 
 They are defined independently. 
\subsection{Contact round surgery of index~$1$} \label{sec:1crdsrg}
  We define contact round surgery of index~$1$. 
 As is mentioned above, a round surgery of index~$1$ is operated 
on two open solid tori in the given $3$-manifold. 
 An open solid torus can be regarded as a tubular neighborhood of a knot. 
 In the case of a contact round surgery, 
it is operated along a transverse link with two components. 
 Let $(M,\xi)$ be a contact $3$-manifold, 
and $\Gamma=\gamma_1\sqcup\gamma_2\subset(M,\xi)$ a transverse link, 
where $\gamma_1$, $\gamma_2$ are two connected components. 

\begin{figure}[htb]
  \centering
  \includegraphics[height=3cm]{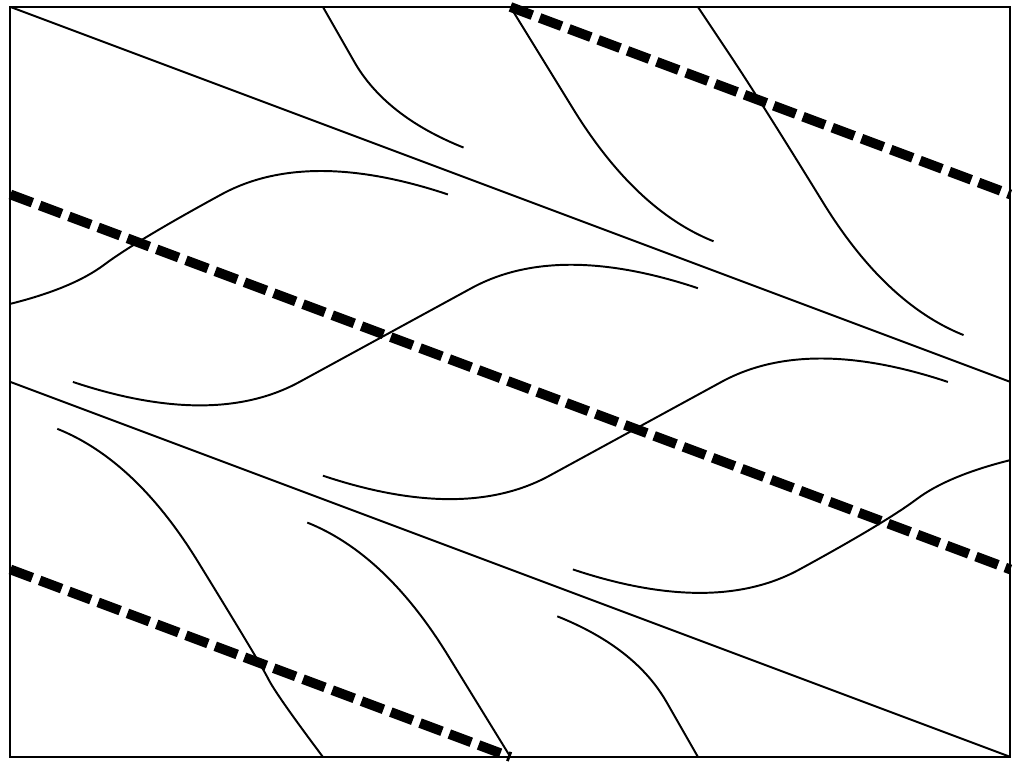}
  \caption{parallel dividing curves}
  \label{fig:pldiv}
\end{figure}

  First, we determine the solid tori to remove. 
 Since each $\gamma_i$ is a transverse knot, 
we can apply the local triviality theorem. 
 According to Theorem~\ref{thm:s1darboux}, there exist contact embeddings 
$\phi_i\colon(S^1\times D(\epsilon_i),\xi_0)\to(M,\xi)$ 
which map $S^1\times\{0\}$ to $\gamma_i$ respectively, 
where $D(\epsilon)$ is a disk with radius $\epsilon$. 
 We may assume $\im\phi_1\cap\im\phi_2=\emptyset$ 
by taking $\epsilon_i$ sufficiently small. 
 We should be careful with the characteristic foliation 
$(S^1\times\rd D(\epsilon_i))_{\xi_0}$ on the boundary tori. 
 It is linear (pre-Lagrangian) with slope $-\epsilon^2$ with respect 
to the meridian $\{\phi=0\}$ and the longitude $\{x=\epsilon,\ y=0\}$ 
of the solid torus. 
 Perturbing the torus slightly, we obtain a convex torus 
with even number of parallel dividing curves 
by Theorem~\ref{thm:convFlex} (see Figure~\ref{fig:pldiv}). 
%
%
 Let $\tilde{\phi}_i\colon(S^1\times\tilde{D}(\epsilon_i),\xi_0)\to(M,\xi)$ 
denote the perturbed ones. 

  Then remove $\inter(\im\tilde{\phi}_1)$ and $\inter(\im\tilde{\phi}_2)$ 
from the given contact $3$-manifold $(M,\xi)$, 
each of which has even number of dividing curves. 

  Next, we construct a contact structure on the thickened torus 
$D^1\times\rd D^2\times S^1=I\times T^2$ 
which is suitable to be glued with 
$M\setminus\{\inter(\im\tilde{\phi}_1)\sqcup\inter(\im\tilde{\phi}_2)\}$. 
 Let $\xi_0=\ker(d\phi-ydx+xdy)=\ker(d\phi+r^2d\theta)$ 
be the standard contact structure on $S^1\times\R^2$. 
 Then the characteristic foliation on the torus 
$T:=\{(\phi,r,\theta)\in S^1\times\mathbb{R}^2\mid r=1\}$ 
is linear with slope $-1$ with respect to the meridian $\{\phi=0\}$ 
and the longitude $\{\theta=0\}$. 
 Perturbing $T$ slightly, we has a convex torus $\tilde{T}$ 
with even number of parallel dividing curves. 
 Since $\tilde{T}$ is convex, 
we have a vertically invariant tubular neighborhood 
$\mathbb{T}\cong T^2\times[-1,1]$ of $\tilde{T}$, 
by using a contact vector field transverse to $\tilde{T}$. 
 Note that the both boundary tori $\rd\mathbb{T}$ 
corresponding to $T^2\times\{-1,1\}$ have the same dividing sets as $\tilde{T}$,
that is, even number of parallel dividing curves non-contractible. 

  Then, according to Proposition~\ref{prop:chfol} 
and Theorem~\ref{thm:convFlex}, 
the thickened torus $(\mathbb{T},\xi_0|_\mathbb{T})$ can be glued to 
$(M\setminus\{\inter(\im\tilde{\phi}_1)\sqcup\inter(\im\tilde{\phi}_2)\},\xi)$ 
according to the dividing curves. 
 Thus we have constructed a contact round surgery of index~$1$ 
along a transverse link with two components. 

  Last of all, we should remark on framings of surgeries. 
 In other words, the choices of coordinates or longitudes of the boundaries 
of the standard tubular neighborhood of transverse knots. 
 The resulting manifold of a round surgery depends relatively on 
both framings of two knots. 
 We can realize any round surgery as a contact surgery above. 
 In fact, the slope of the characteristic foliation on the boundary 
of the standard tubular neighborhood of a transverse knot 
can be taken arbitrary close to $0$ for any framing, 
by taking the tubular neighborhood sufficiently close to the transverse knot. 
 Then we can take any relative framings for surgeries. 

\subsection{Contact round surgery of index~$2$} 
\label{sec:2crdsrg}
  We define contact round surgery of index~$2$. 
 First of all, we should remark that this surgery is not always defined. 
 As is mentioned in the last part of Section~\ref{sec:rd-handle}, 
a round surgery of index~$2$ is operated 
along a $2$-dimensional torus embedded in the given $3$-manifold. 
 An open thickened torus as a tubular neighborhood of the embedded torus 
is removed. 
 And then two solid tori are reglued. 
 In the case of a contact round surgery of a contact $3$-manifold $(M,\xi)$, 
it is operated along an embedded torus $T\subset(M,\xi)$.  
 We impose this torus the condition that
\begin{equation} \label{eq:releul}
  \langle e(\xi),[T]\rangle=0, 
\end{equation}
where $e(\xi)\in H^2(M;\Z)$ is the Euler class of the contact structure, 
and $[T]\in H_2(M;\Z)$ is the fundamental class. 
 This condition is translated to Equation~\eqref{eq:eulcond} 
in terms of convex surface theory if $T$ is a convex torus. 
 In the first part of this subsection, we discuss on this condition. 
 And then, under this condition, we define contact round surgery of index~$2$. 

\subsubsection{The condition on the embedded tori} \label{sec:conditorus}
  Now, we discuss the meaning of the condition 
on an embedded torus $T\subset(M,\xi)$ 
where a contact round surgery is operated. 
 The condition is required in the construction of a contact structure 
on the surgered manifold. 
 In the following construction, 
we remove a tubular neighborhood of $T\subset(M,\xi)$, 
and reglue two contact solid tori. 
 In order to do that, we need contact structures on a solid torus 
which have the same characteristic foliation (or dividing set) 
as $T\subset(M,\xi)$. 
 In other words, we need to extend the contact structure 
determined by the characteristic foliation to the whole solid torus. 
 The Euler class $e(\xi)\in H^2(M;\Z)$ of $\xi$ 
evaluated with $[T]\in H_2(M;Z)$ is an obstruction 
to the extension as a plane field by the obstruction theory 
(see for example \cite{bredon}). 
 That is, Condition~\eqref{eq:releul} guarantees 
the extension as a plane field. 
 Actually, contact structure along $T$ extends as a contact structure 
in this case. 
 It is proved in the construction below. 

  Condition~\eqref{eq:releul} for a $2$-dimensional torus $T$ 
embedded in a contact $3$-manifold $(M,\xi)$ 
is natural in some cases. 
 We have the following two examples. 

%
%
\begin{expl} \label{expl:tight}
  If $\xi$ is a tight contact structure on a $3$-manifold $M$, 
any embedded torus $T\subset(M,\xi)$ satisfies Condition~\eqref{eq:releul}. 
 In fact, for tight contact structures, the following property is known: 
%
%
\begin{thrm}[Eliashberg, \cite{eliash20years}]
  Let $(M,\xi)$ be a tight contact $3$-manifold. 
 For any closed orientable surface $\Sigma$ embedded in $(M,\xi)$, 
the following inequality holds\textup{:} 
\begin{equation*}
  |\langle e(\xi),[\Sigma]\rangle|\le
  \begin{cases}
    \ 0 & \text{if}\ \Sigma=S^2, \\
    \ -\chi(\Sigma) & \text{otherwise}, 
  \end{cases}
\end{equation*}
where $\chi(\Sigma)$ is the Euler characteristic of $\Sigma$. 
\end{thrm} 
\noindent
 In the case that we consider now, 
$\Sigma$ is a torus $T$ and then $\chi(T)=0$. 
 Therefore, $\releu{T}=0$. \qed
\end{expl}
%
%
\begin{expl} \label{expl:euler}
  If $M$ is a closed orientable $3$-manifold, 
any embedded torus $T\subset(M,\xi)$ which separates the manifold $M$ 
satisfies Condition~\eqref{eq:releul}. 
 In fact, 
$\releu T$ is an obstruction 
to extending the contact plane field $\xi|_T$ along $T$ 
to the 
manifold bounded by $T$ as a plane field. 
 In this case, $\xi$ itself is the extension. 
 Then $\releu T$ should vanish. 
\qed
\end{expl} 

 The case when the contact structure $\xi$ on a closed $3$-manifold $M$ 
is overtwisted 
and the embedded torus $T\subset M$ does not separate 
the underlying manifold is left. 
 In that case, we can not define contact round surgery of index~$2$ 
operated along the torus, which has $\releu T\ne0$. 
 However, after modifying the contact structure $\xi$ suitably 
to $\xi'$ on $M$ with $\releu[\xi']T=0$, 
we can operate a contact round surgery of index~$2$ along the same $T$. 
 we observe these operation in Section~\ref{sec:concl}. 

\subsubsection{Definition of a contact round surgery 
  of index~$2$}
\label{sec:defind2}
  Now, we define contact round surgery of index~$2$ 
under Condition~\eqref{eq:releul}. 
 Let $(M,\xi)$ be a contact $3$-manifold, 
and $T\subset(M,\xi)$ an embedded $2$-dimensional torus 
along which the round surgery is operated. 
 Assume that $T$ satisfies Condition~\eqref{eq:releul}. 
 We may also assume that $T\subset(M,\xi)$ is convex 
due to Theorem~\ref{thm:convGen}. 
 In addition, we take a meridian $\hat{\mu}\subset T$, 
or $\hat{\mu}\in H_1(T,\Z)$. 
 It corresponds to a meridian 
$\rd D^2\times\{\ast\}\times\{\ast\}\subset\rd D^2\times D^1\times S^1$ 
of the thickened torus to be removed. 
 In other words, it will be the meridians 
$\rd D^2\times \rd D^1\times\{\ast\}\subset D^2\times\rd D^1\times S^1$ 
of the solid tori to be reglued. 
 In the following, we call it the \emph{surgery meridian}. 

  First, we perturb the embedded torus $T$ to a suitable position. 
 We isotope $T$ so that the dividing set $\Gamma_T$  on $T$ 
consists of homotopically non-trivial parallel curves. 
 In order to do that, we should check possible dividing sets 
on a convex $2$-dimensional torus under Condition~\eqref{eq:releul}. 
 According to the properties of dividing sets, 
a dividing set on a convex $2$-dimensional torus consists 
of even-number of parallel homotopically non-trivial curves 
and some homotopically trivial curves who never intersect each other. 
 Further, on the Euler characteristics of the domains $U_\pm\subset T^2$ 
divided by dividing sets, 
there exists the following constraint: 
\begin{equation} \label{eq:eulcond}
  \chi(U_+)=\chi(U_-)=0. 
\end{equation}
 In fact, we have the following formulas 
for a convex closed surface $\Sigma$ in a contact $3$-manifold $(N,\tilde\xi)$ 
with positive and negative regions $R_\pm$ (see \cite{hondaI}): 
\begin{equation*} \label{eq:formulae}
  \chi(\Sigma)=\chi(R_+)+\chi(R_-),\qquad 
  \releu[\tilde\xi]{\Sigma}=\chi(R_+)-\chi(R_-). 
\end{equation*}
 In this case, $\chi(\Sigma)=\chi(T)=0$ since $T$ is a torus, 
and $\releu T=0$ from the assumption. 
 Then we have $\chi(U_+)=\chi(U_-)=0$. 

  Then we consider removing homotopically trivial dividing curves 
in what follows. 
 The method to be applied is the bypass attachment. 
 Now that the dividing set $\Gamma_T$ on $T$ 
has a homotopically trivial dividing curves, 
on the transversely invariant tubular neighborhood $U\subset M$ of $T$, 
the contact structure $\xi$ is overtwisted 
from the Giroux criterion (Proposition~\ref{prop:girouxcrit}). 
 Further, $\xi$ is overtwisted on $U\setminus T$. 
 Then, by Proposition~\ref{prop:huang}, we can find any embedded bypass we need.
 Isotoping $T$ along the bypass, the dividing set $\Gamma_T$ is modified 
in the same manner as the bypass attachment 
(See Lemma~\ref{lem:bpatt}, Figure~\ref{fig:bypatt}). 
 Therefore, we have only to follow such modifications of the dividing set 
in order to cancel homotopically trivial dividing curves. 

\begin{figure}[htb]
  \centering
  \includegraphics[height=4cm]{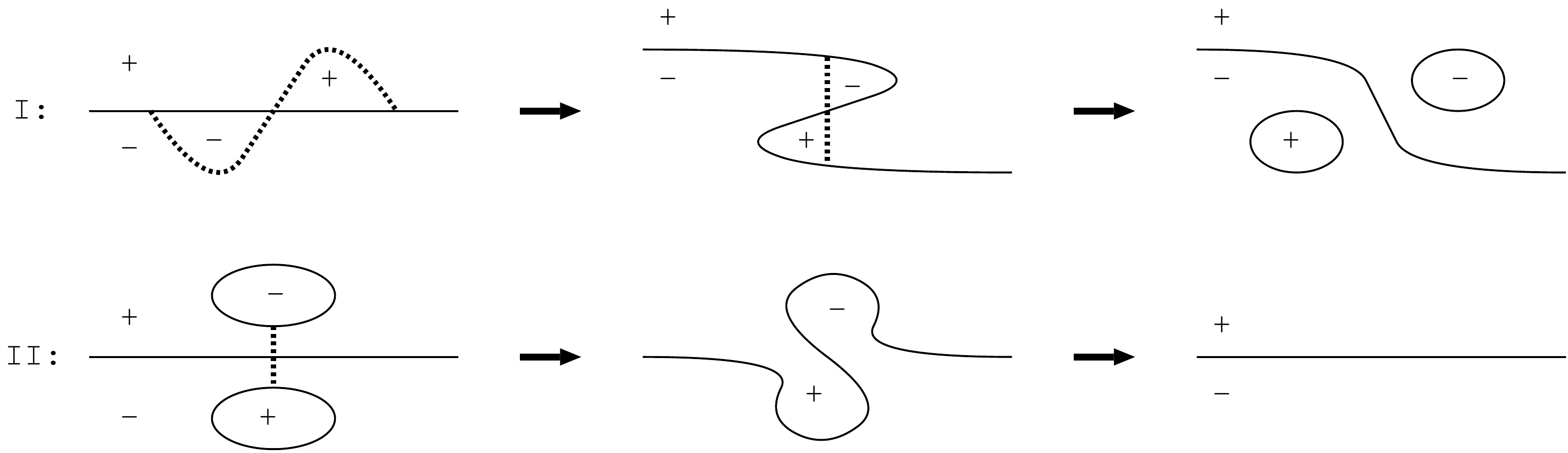}
  \caption{creating and canceling independent contractible dividing curves. }
  \label{fig:ctrdiv}
\end{figure}

\begin{figure}[htb]
  \centering
\includegraphics[width=14cm]{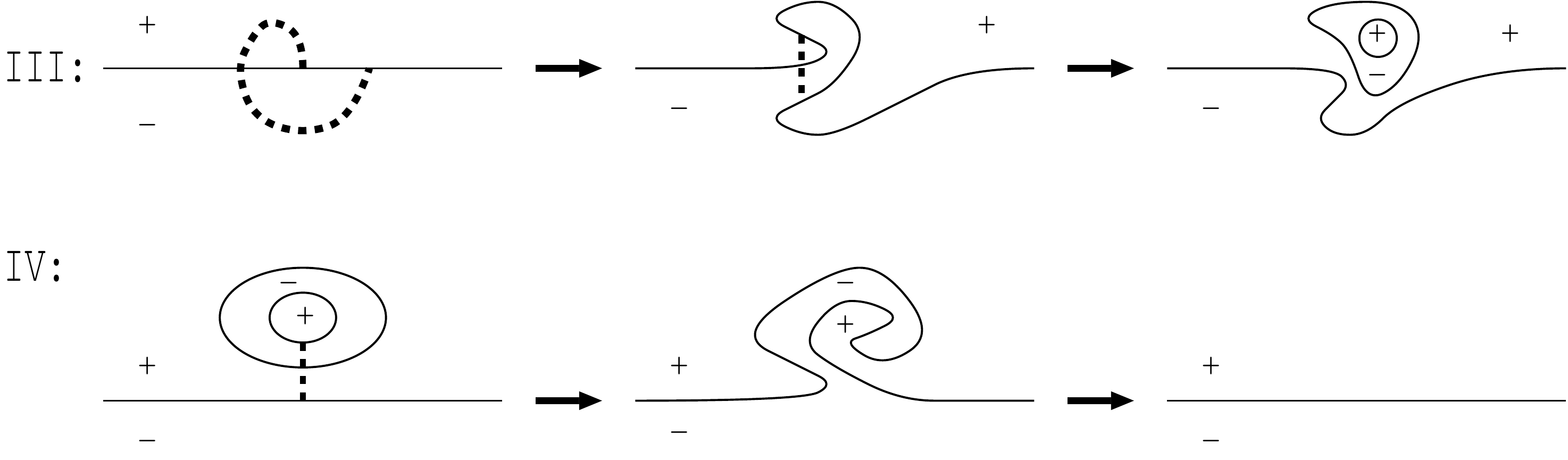}  
  \caption{nested pair of dividing curves}
  \label{fig:nstddiv}
\end{figure}

  We introduce basic operations creating or canceling 
a pair of homotopically trivial dividing curves. 
 By attaching a bypass, we can create or cancel a pair of homotopically trivial 
dividing curves. 
 See Figure~\ref{fig:ctrdiv} for an independent pair 
and Figure~\ref{fig:nstddiv} for a nested pair. 
%
%
%
%

\begin{figure}[htb]
  \centering
  \includegraphics[width=14cm]{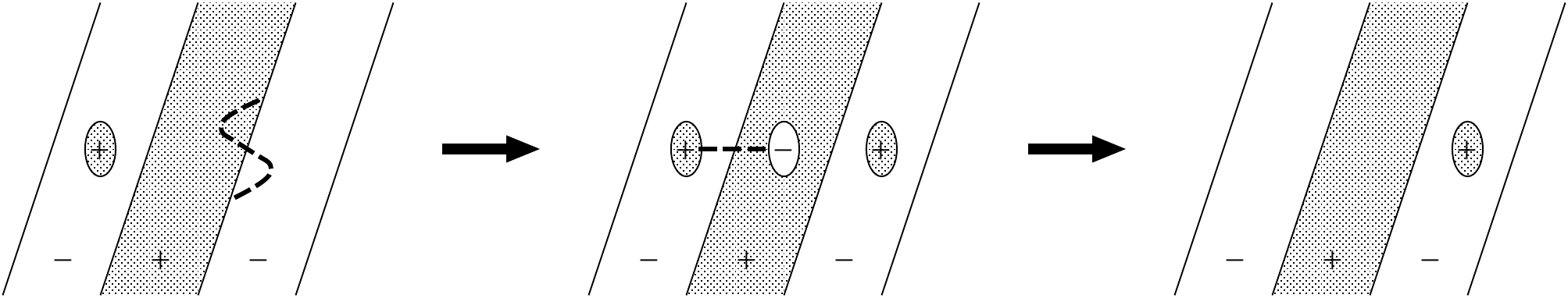}
  \caption{move to other domain}
  \label{fig:mv2odom}
\end{figure}

  By using Operations~\texttt{I}, \texttt{II}, \texttt{III}, and \texttt{IV}, 
any possible dividing set is reduced 
to homotopically non-trivial dividing curves. 
 First, we apply Operation \texttt{IV} 
to cancel nested pair of homotopically trivial dividing curves. 
 Then the rest of homotopically trivial dividing curves 
are single independent ones. 
 Such a single homotopically trivial dividing curve is moved to other domain 
by the combination of Operations~\texttt{I} and \texttt{II} 
(see Figure~\ref{fig:mv2odom}). 
%
%
 All homotopically trivial dividing curves which bound positive domain 
are gathered in one negative strip , 
and all those which bound negative domain are gathered in a positive strip 
next to the positive strip above. 
 From Condition~\eqref{eq:eulcond}, 
the numbers of homotopically trivial dividing curves 
bounding positive or negative domains must be the same. 
 Therefore, all of them are canceled by Operation~\texttt{II}. 
 Then the given contact solid torus with convex boundary is reduced 
to the contact solid torus without homotopically trivial dividing curves 
on the boundary. 

  Now, we operate a round surgery of index~$2$ 
along the isotoped convex torus $T$ 
with no homotopically trivial dividing curve. 
 Recall that the dividing set $\Gamma_T$ on $T$ consists of even number of 
parallel dividing curves. 

 First, we determine the thickened torus to remove. 
 Since $T\subset(M,\xi)$ is convex, 
there exists a contact vector field $X$ transverse to $T$ 
on a neighborhood of $T\subset M$. 
 By using the flow of the vector field $X$, an open tubular neighborhood 
$\phi\colon T\times(-\epsilon,\epsilon)\to(M,\xi)$ is constructed, 
where $\phi$ maps $T\times\{0\}$ to $T\subset M$ identically. 
 Remove $\im\phi\subset M$ from the given contact $3$-manifold $(M,\xi)$. 
 We should be careful with the characteristic foliation on the boundary 
of $M\setminus\im\phi$. 
 It is diffeomorphic to the two copies 
of the characteristic foliation $T_\xi$ because the vector field $X$ is contact.

  Next, we need two contact solid tori to reglue. 
 Recall that we have a surgery meridian $\hat{\mu}$ of $T$. 
 For one of the contact solid torus it should be the meridian. 
 On the other hand, for the other one, $-\hat{\mu}$ 
should be the meridian. 
 In other words, the characteristic foliations 
on the boundary tori are the mirror image of each other. 
 The model is constructed as follows. 
 Let $\zeta$ be an overtwisted contact structure 
$\{(\cos r^2)d\phi+r(\sin r^2)d\theta=0\}$ on $S^1\times\mathbb{R}^2$ 
with cylindrical coordinates $(\phi,r,\theta)$ 
as in Section~\ref{sec:1crdsrg}. 
 As is seen above, the characteristic foliation $(S^1\times \rd D(\rho))_\zeta$ 
is a linear foliation with slope $-\tan\rho^2$. 
 By perturbing it, we have even number of parallel dividing curves 
by Theorem~\ref{thm:convFlex} (see Figure~\ref{fig:pldiv}). 
 For $0<\rho<\pi$, we obtain any non-zero slope. 
 Even in the case of meridional dividing curves, 
we obtain them by perturbing $S^1\times\rd D(\sqrt{\pi})$. 
 Thus we obtain any required homotopically non-trivial dividing curves 
by perturbing $S^1\times\rd D(\rho)\subset (S^1\times\R^2,\zeta)$, 
$0<\rho\le\sqrt{\pi}$. 
 Then, if the slope is not $0$, the model is considered 
as a part of a contact $3$-manifold $(S^1\times S^2,\zeta_0)$ 
which is a union of the two copies of $(S^1\times D(\sqrt{\pi/2}),\zeta)$. 
 Even if the slope is $0$, it is considered 
as a part of $(S^1\times S^2,\zeta_1)$ 
which is from the two copies of $(S^1\times D(\sqrt{\pi}),\zeta)$. 
 Therefore, the closure of the complement of the obtained solid torus 
is the other required one. 

  We glue the obtained contact solid tori  to $(M\setminus\im\phi,\zeta)$ 
as follows. 
 From the construction, they have the same dividing curves. 
 By Theorem~\ref{thm:convFlex}, we may regard that they have 
the same characteristic foliations. 
 Then, by Proposition~\ref{prop:chfol}, we can glue them. 

  Last of all, we remark that  
we can define contact round surgery for any framing, 
or surgery meridian,  of the embedded torus 
because we can arrange the slope of the dividing curves 
of the model of contact solid torus for that. 

\section{Lutz twist and Giroux torsion} \label{sec:LtwGtor}
  There exist two important notions on contact structures on $3$-manifolds, 
the so called Lutz twist and Giroux torsion. 
 First, we recall the definitions of them. 
 Then we show that the Lutz twist is realized by contact round surgeries. 
 Further, we discuss some relation to the Giroux torsion. 

\subsection{Definitions} \label{sec:dfnlztw}
¡¡¡¡In this subsection, we review the definitions
of the Lutz twists and the Giroux torsion. 

\subsubsection{Lutz twist along $S^1$}
\begin{figure}[htb]
  \centering
  \def\svgwidth{10.5cm}
  {\small 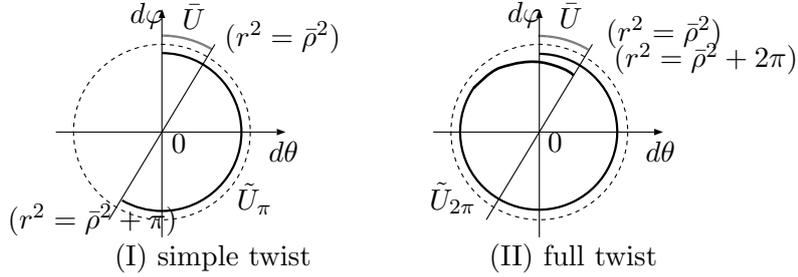}
  \caption{\vtop{\hbox{Lutz twist along $S^1$} 
    \hbox{(situations of covectors or contact form $\zeta$)}}}
  \label{fig:lztwdef}
\end{figure}

  Let us begin by reviewing the definition of the original Lutz twist. 
 It is an operation modifying a contact structure along a transverse knot. 
 Let $\Gamma$ be a transverse knot in a contact $3$-manifold $(M,\xi)$. 
 There exists a tubular neighborhood of $\Gamma$ which is contactomorphic 
to $(S^1\times\dsc\rho,\xi_0)=:U$ for some small radius $\rho>0$,
by Theorem~\ref{thm:s1darboux}. 
 In order to define the Lutz twist, 
we need the standard overtwisted contact structure. 
 Let $\zeta$ be an overtwisted contact structure on $S^1\times\mathbb{R}^2$ 
defined as 
\begin{equation*}
  \zeta:=\ker\{(\cos r^2)d\phi+(\sin r^2)d\theta\},
\end{equation*}
where $(\phi,r,\theta)\in S^1\times\mathbb{R}^2$ 
are the cylindrical coordinates. 
 Note that $(\stm\R^2,\xi_0)$ is isotopic to $(\stm\inter\dsc{\pi/2},\zeta)$. 
 By this correspondence, the tubular neighborhood $U=(\stm\dsc{\rho},\xi_0)$ 
is identified with $(\stm\dsc{\bar\rho},\zeta)=:\bar{U}$ for $\bar\rho>0$ 
satisfying $\rho^2=\tan\bar\rho^2$, $0<\bar\rho^2<\pi/2$. 
 The \emph{simple Lutz twist}\/ is the operation replacing 
$U\cong\bar{U}=(\stm\dsc{\bar\rho},\zeta)$ 
with $(S^1\times\dsc{\sqrt{\bar\rho^2+\pi}},\zeta)=:\tilde{U}_\pi$ 
(see Figure~\ref{fig:lztwdef}-(I)). 
%
%
 In fact, since the characteristic foliations 
$(\rd\bar{V})_\zeta=(S^1\times\rd D(\bar\rho))_{\zeta}$ 
and $(\rd\tilde{V}_\pi)_\zeta=(S^1\times\rd D(\sqrt{\bar\rho^2+\pi}))_\zeta$ 
are the same, 
these two contact solid tori can be replaced. 
 Note that the contact plane of the second one twists 
a half time, or $\pi$, more than the first one along a radial $r$-axis.  
 Similarly, the \emph{full Lutz twist}\/ is defined as the operation replacing 
$V=(S^1\times\dsc\rho,\xi_0)$ 
with $\tilde{V}_{2\pi}:=(S^1\times\dsc{\sqrt{\bar\rho+2\pi}},\zeta)$
(see Figure~\ref{fig:lztwdef}-(II)). 
 In this case, the contact plane twists one time, or $2\pi$, more. 
 We should remark here that a simple Lutz twist does change 
the homotopy class of the contact structure as plane fields, 
a full Lutz twists does not though. 
 It is clear that a Lutz twist make a contact structure be overtwisted. 

\subsubsection{Lutz twist along $T^2$} \label{sec:lztwalgt2}
\begin{figure}[htb]
  \centering
  \def\svgwidth{10.5cm}
  {\small 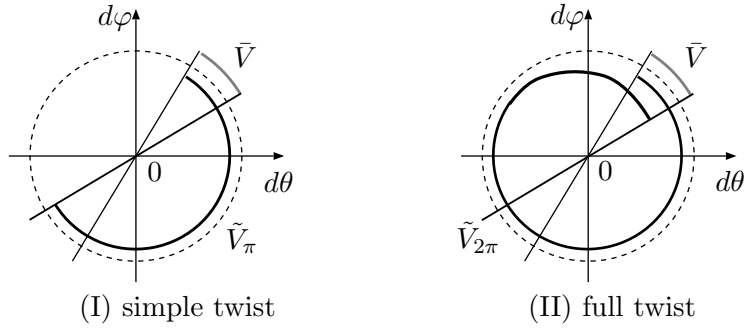}
  \caption{Lutz twist along a torus}
  \label{fig:glztwdef}
\end{figure}

  Next, we introduce the Lutz twist along a certain torus. 
 Let $T$ be a torus embedded in a contact $3$-manifold $(M,\xi)$. 
 Assume that it is pre-Lagrangian. 
 An embedded torus $T\subset(M,\xi)$ is said to be \emph{pre-Lagrangian}\/ 
if the characteristic foliation $T_\xi$ is linear with closed leaves. 
 In other words, $T$ is foliated by Legendrian circles. 
 Then it is well known (see \cite{geitext}, \cite{giroux99}) 
that it has a tubular neighborhood which is contactomorphic 
to a tubular neighborhood of the $2$-dimensional torus 
$\{(\phi,r,\theta)\in S^1\times\R^2\mid r=\rho\}\subset(S^1\times\R^2,\xi_0)$, 
for some radius $\rho>0$. 
 We may regard the tubular neighborhood as 
\begin{align*}
  &V:=(\{(\phi,r,\theta)\mid\rho-\delta_1<r<\rho+\delta_2\},\xi_0)\\
  \cong&\bar{V}:=(\{(\phi,r,\theta)\mid
  \bar{\rho}^2-\delta<r^2<\bar{\rho}^2+\delta\},\zeta),
\end{align*}
where $\rho^2=\tan\bar\rho^2$, $(\rho-\delta_1)^2=\tan(\bar\rho^2-\delta)$,
$(\rho+\delta_2)^2=\tan(\bar\rho^2+\delta)$. 
 The \emph{simple Lutz twist}\/ (or \emph{$\pi$-Lutz twist}\/) 
\emph{along}\/ $T\subset(M,\xi)$ is defined 
as the operation replacing the tubular neighborhood $\bar{V}$ of $T$ with 
\begin{equation*}
  \tilde{V}_\pi:=
  (\{(\phi,r,\theta)\mid\bar\rho^2-\delta<r^2<\bar\rho^2+\pi+\delta\}, \zeta), 
\end{equation*}
(see Figure~\ref{fig:glztwdef}). 
%
%
 We can define the \emph{full Lutz twist}\/ (or \emph{$2\pi$-Lutz twist}\/) 
\emph{along}\/ $T$
by using $\tilde{V}_{2\pi}:=\{(\phi,r,\theta)\mid
\bar\rho^2-\delta<r^2<\bar\rho^2+2\pi+\delta\}$. 

  The Lutz twist along a pre-Lagrangian torus is considered 
as a generalization of that along a transverse knot.
 In fact, for a transverse knot $\Gamma\subset(M,\xi)$, 
there exists a tubular neighborhood $U\subset(M,\xi)$ of $\Gamma$
which is contactomorphic to a tubular neighborhood
of $\stm\{0\}\subset(\stm\R^2,\xi_0)$.
 Then, in the tubular neighborhood $U$,
we have a pre-Lagrangian torus $T\subset U$ corresponding
to some torus $\{(\phi,r,\theta)\mid r=\rho\}\subset(\stm\R^2,\xi_0)$.
 The Lutz twist along $T$ is equivalent to that along $\Gamma$.

 We should mention that the Lutz twist along a pre-Lagrangian torus
may not create any overtwisted disc. 
 However, the Lutz twist along a pre-Lagrangian torus creates 
the following important thing. 

\subsubsection{Giroux torsion}\label{sec:Gtor}
  We review the definition of the Giroux torsion. 
 A contact manifold $(M,\xi)$ is said 
to have the \emph{Giroux torsion at least} $n\in\mathbb{N}$ 
if there exists a contact embedding 
$f_n\colon(T^2\times I, \tilde\zeta_n)\to(M,\xi)$, 
where $\tilde\zeta_n$ is a contact structure on $T^2\times I$ with coordinates 
$(\phi,r,\theta)\in S^1\times I\times S^1\subset S^1\times\mathbb{R}^2$ 
defined as 
\begin{equation*}
  \tilde\zeta_n:=\ker\{\cos(2n\pi r)d\theta+\sin(2n\pi r)d\phi\}.   
\end{equation*}
 The supremum of these numbers for all such embeddings 
to the contact manifold $(M,\xi)$ 
is called the \emph{Giroux torsion}\/ of $(M,\xi)$. 
 Let $\operatorname{Tor}(M,\xi)$ denotes it. 
 If there exists no such embedding, 
$(M,\xi)$ is said to have Giroux torsion $0$. 
 The definition can be extended for half-integers $n=m/2$, $m\in\mathbb{N}$. 
 Especially, for $n=1/2$, 
we call the contact embedding 
$f_{1/2}\colon(T^2\times I, \tilde\zeta_n)\to(M,\xi)$ 
a \emph{half Giroux torsion unit}. 

  The Giroux torsion is an invariant of contact $3$-manifold 
introduced by Giroux~\cite{girouxTor}, 
which explicitly appears in the classification of tight contact structures 
on the $3$-dimensional torus $T^3$ 
due to Giroux and Kanda~\cite{kanda} independently. 
 It is proved that a closed contact manifold 
is not strongly symplectically fillable 
if it has the Giroux torsion greater than $0$ (see~\cite{gay06}). 

  The Giroux torsion is closely related to the Lutz twist. 
 The Lutz twist along a torus make a Giroux torsion unit. 
 In fact, comparing the contact structures $\tilde{\zeta}_n$ and $\zeta$, 
we obtain that the substitute thickened torus in the Lutz twist includes 
a half Giroux torsion unit. 

\subsection{Realization by round surgeries}\label{sec:cstrlztw}
  Now, we show that the Lutz twist 
is realized by a certain $4$-tuple of contact round surgeries. 
 More precisely, we construct two pairs of contact round surgeries
of index~$1$ and index~$2$.
 By such contact round surgeries,
we will realize the Lutz twist along a pre-Lagrangian torus.
 As is mentioned in Subsubsection~\ref{sec:lztwalgt2},
this realizes the Lutz twist along a transverse knot as well. 
 The claim in this subsection is the following. 
%
%
\begin{thrm} \label{thm:Lz}
  A simple Lutz twist is realized by a certain ordered $4$-tuple 
of a contact round surgeries of index~$1$ and those of index~$2$. 
\end{thrm} 

  First of all, we confirm the starting situations locally.
 Let $T\subset(M,\xi)$ be a pre-Lagrangian torus
in a contact $3$-manifold $(M,\xi)$ 
where we will operate the Lutz twist. 
 Along the pre-Lagrangian torus $T$, 
there exist the standard tubular neighborhood $U\subset(M,\xi)$ 
which is contactomorphic to a tubular neighborhood 
of the torus $\{(\phi,r,\theta)\mid r=\rho\}\subset
(S^1\times\R^2,\zeta=\ker\{(\cos r^2)d\phi+(\sin r^2)d\theta\})$ 
for some $\rho>0$. 
 We may assume that the tubular neighborhood is 
$U=(\{(\phi,r,\theta)\mid\pi/4-\epsilon<r^2<\pi/4+\epsilon\},\zeta)$ 
by taking longitude and meridian suitably. 
 In the following, we discuss by using this local model.  

  We operate a contact round surgery of index~$2$ along $T\subset(U,\zeta)$ 
first. 
 Since the torus clearly satisfies Condition~\eqref{eq:releul}, 
we can operate a contact round surgery of index~$2$ along $T$. 
 In the local model $(U,\zeta)$, we determine the surgery meridian as follows. 
 Let $\mu\in H_1(T;\Z)$, or $\mu\subset T$, be the meridian of $T$ 
corresponding to $\{\phi=0\}$ 
and $\lambda\in H_1(T;\Z)$, or $\lambda\subset T$, the longitude of $T$ 
corresponding to $\{\theta=0\}$. 
 We take $\hat{\mu}:=\lambda$ as the surgery meridian. 
 Note that, by the framing $(\mu,\lambda)$ of $T$, the characteristic foliation 
$T_\zeta$ on $T$ is represented by $\lambda-\mu$. 

\begin{figure}[htb]
  \centering
  \def\svgwidth{11.1cm}
  {\small 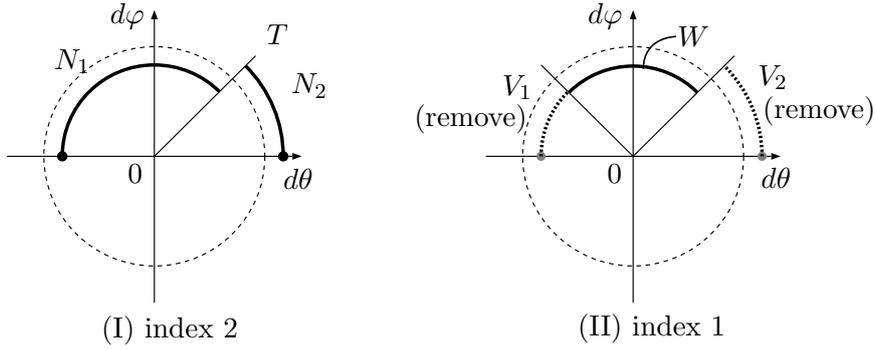}
  \caption{\vtop{\hbox{Contact round surgeries}
      \hbox{(situations of covectors or contact forms)}}}
  \label{fig:lztwctsg}
\end{figure}

  Now, we operate a contact round surgery of index~$2$ along $T\subset(U,\zeta)$
with surgery meridian $\hat{\mu}=\lambda$. 
 Cutting $(U,\zeta)$ open along $T$, we reglue two contact solid tori. 
 These tori are prepared according to the surgery meridian $\hat{\mu}$ 
and the characteristic foliation $T_\zeta$. 
 Let $N_1\cong S^1\times D^2$ denote a solid torus 
whose boundary has the same orientation as $T$, 
and $N_2$ denote a solid torus 
whose boundary has the opposite orientation to $T$. 
 For $N_1$, the meridian $\mu_1$ of the boundary torus 
$\rd N_1\cong S^1\times\rd D^2$ corresponds to $\hat{\mu}=\lambda$. 
 As its longitude, we take $\lambda_1\subset\rd N_1$, 
or $\lambda_1\in H_1(\rd N_1;\Z)$, that corresponds to $-\mu\in H_1(T;\Z)$, 
so that the orientation of $\rd N_1$ is the same as $T$. 
 Then the characteristic foliation $T_\zeta$ corresponds 
to the foliation on $\rd N_1$ represented by $\mu_1+\lambda_1$. 
 Therefore, $N_1$ should be isotopic to 
$(S^1\times D(\sqrt{3\pi/4}),\zeta)$ (see Figure~\ref{fig:lztwctsg}-(I)). 
 Similarly, for $N_2$, the meridian $\mu_2$ of $\rd N_2$ 
corresponds to $-\hat{\mu}=-\lambda$ 
and a longitude $\lambda_2$ corresponds to $-\mu$, 
so that the orientation of $\rd N_2$ is opposite to $T$. 
 Then the characteristic foliation $T_\zeta$ corresponds 
to the foliation on $\rd N_2$ represented by $-\mu_2+\lambda_2$. 
 Therefore, $N_2$ should be isotopic to 
$(S^1\times D(\sqrt{\pi/4}),\zeta)$ 
(see Figure~\ref{fig:lztwctsg}-(I)). 
 Gluing these two solid tori $N_1$ and $N_2$ to $(U\setminus T,\zeta)$, 
we obtain a new contact manifold. 

  Next, we operate a contact round surgery of index~$1$. 
 Recall that in the previous contact round surgery of index~$2$, 
we glue two contact solid tori $N_1$ and $N_2$. 
 Let $\gamma_1\subset N_1$ and $\gamma_2\subset N_2$ 
be transverse knots corresponding 
to $S^1\times\{0\}\subset(S^1\times\dsc{\sqrt{3\pi/4}},\zeta)=N_1$, 
$(S^1\times\dsc{\sqrt{\pi/4}},\zeta)=N_2$, respectively. 
 We operate a contact round surgery of index~$1$ 
along the transverse link $\gamma_1\sqcup\gamma_2$. 
 In other words, removing tubular neighborhoods of $\gamma_1$ and $\gamma_2$, 
we glue two boundary tori together. 
 In order to determine the framing of surgery, 
we determine tubular neighborhoods of $\gamma_1$ and $\gamma_2$. 
 As the tubular neighborhoods, we take 
$V_1:=(\stm\dsc{\sqrt{\pi/4}},\zeta)\subset N_1$ for $\gamma_1$ 
and $V_2:=(\stm\dsc{\sqrt{\pi/4}},\zeta)\subset N_2$ for $\gamma_2$ 
(see Figure~\ref{fig:lztwctsg}-(II)). 
 Then removing $\inter V_1$ and $\inter V_2$, 
we glue two boundary tori $\rd V_1$ and $\rd V_2$ 
so that their meridians and characteristic foliations agree. 
 Thus we obtain a new contact manifold. 

%
%

  As a result, these two contact round surgeries amount to the following. 
 In total, cutting the contact manifold $(U,\zeta)$ open 
along the pre-Lagrangian torus $T\subset(U,\zeta)$, 
we glued the boundary tori to the contact thickened torus 
$\{(\phi,r,\theta)\mid 0\le r^2\le\pi/2\}\cong T^2\times I$ 
from the both sides. 
 In fact, in the first contact round surgery of index~$2$, 
cutting $(U,\zeta)$ open along $T$, we glued two solid tori $N_1$ and $N_2$. 
 Then, in the second contact round surgery of index~$1$, 
removing  
\begin{align*}
  \inter V_1&=(\stm\inter\dsc{\sqrt{\pi/4}},\zeta)
              \subset(\stm\dsc{\sqrt{3\pi/4}},\zeta)= N_1, \\
  \inter V_2&=(\stm\inter\dsc{\sqrt{\pi/4}},\zeta) 
              \subset(\stm\dsc{\sqrt{\pi/4}},\zeta)=N_2, 
\end{align*}
we glued $\rd V_1$ and $\rd V_2$ together. 
 As a result, the contact thickened torus 
\begin{equation*}
  W_1:=N_1\setminus\inter V_1
  =(\{(\phi,r,\theta)\mid\pi/4\le r^2\le3\pi/4\},\zeta)
\end{equation*}
is left in the resulting manifold. 
 It is isotopic to $(\{(\phi,r,\theta)\mid0\le r^2\le\pi/2\},\zeta)$ 
(see Figure~\ref{fig:lztwctsg}-(II)). 
 We should remark that, at this moment, 
the underlying manifold has been modified. 

  We operate the same pair of contact round surgeries again 
along the pre-Lagrangian torus $\tilde{T}=\rd V_1=\rd V_2$, 
where the first pair of round surgeries finished. 
 Then we obtain another contact thickened torus 
$W_2:=(\{(\phi,r,\theta)\mid0\le r^2\le\pi/2\},\zeta)$ 
just next to the previous $W_1$. 

  These two pairs of contact round surgeries amount to the simple Lutz twist 
along $T\subset(U,\zeta)$. 
 In fact, the contact thickened tori $W_1$ and $W_2$ are glued 
so that their meridians and characteristic foliations agree. 
 Then the combined contact thickened torus $W_1\cup W_2$ is isotopic to 
$(\{(\phi,r,\theta)\mid0\le r^2\le\pi\},\zeta)$. 
 We should remark that the underlying manifold has recovered 
to the original shape 
because the boundary tori $\{r^2=0\}$ and $\{r^2=\pi\}$ 
have the same characteristic foliation. 
 This implies that the total operation is nothing but the simple Lutz twist. 

  Thus, Theorem~\ref{thm:Lz} has been proved. 

\section{Conclusion} \label{sec:concl}
  Theorem~\ref{thma} and Theorem~\ref{thmb} are proved in this section. 
 In the first subsection, we show that any closed orientable $3$-manifold admits
a contact structure, by using round surgeries. 
 In the second subsection, we show that any contact structure 
on any closed orientable $3$-manifold is obtained 
from the standard contact $3$-sphere by contact round surgeries.

\subsection{Construction on any manifold} \label{sec:prThma}
  We give a proof of Theorem~\ref{thma}, 
that is, an alternative proof of the theorem 
firstly proved by Martinet~\cite{martinet}.
 In this paper, it is proved by using round surgeries. 
  We use the method as follows. 
 According to Theorem~\ref{thm:r-surg}, 
any closed orientable $3$-manifold is constructed 
from a $3$-dimensional sphere by a sequence of round surgeries 
of index~$1$ and index~$2$. 
 What is to be proved is that the given sequence of round surgeries 
can be operated as contact round surgeries. 
 First, we show that each round surgery of index~$1$ can be operated 
as a contact round surgery. 
 On the other hand, some round surgeries of index~$2$ can not be operated 
as contact round surgeries directly.  
 Recall that a contact round surgery of index~$2$ is defined 
along a torus $T^2$ embedded into a contact $3$-manifold $(M,\xi)$ 
which satisfies Condition~\eqref{eq:releul}: $\releu{T^2}=0$. 
 We should discuss the case 
when the round surgery that we would like to operate as a contact round surgery 
is operated along a torus $T\subset M$ with $\releu T\ne 0$. 
 In that case, we modify the contact structure $\xi$ so that it satisfies 
$\releu T=0$. 
 At the end of this subsection, 
we discuss the existence of a sequence of contact round surgeries 
without modifying the contact structures. 

  First, we deal with round surgeries of index~$1$. 
 A surgery of this kind is operated along a link with two components. 
 According to Proposition~\ref{prop:trsvApp}, any curve is approximated 
by a positively transverse curve. 
 Then any round surgery of index~$1$ is operated as a contact round surgery. 

  The argument on round surgeries of index~$2$ is much more complicated. 
 A round surgery of index~$2$ is operated 
along an embedded $2$-dimensional torus. 
 A contact round surgery of index~$2$ is operated along a torus $T$ 
in a contact $3$-manifold $(M,\xi)$ with Condition~\eqref{eq:releul}: 
$\releu T=0$. 
 By Theorem~\ref{thm:convGen}, the torus $T$ is approximated by a convex torus. 
 Since we deal with a closed contact $3$-manifold $(M,\xi)$, 
Condition~\eqref{eq:releul} is satisfied if $\xi$ is tight 
or $T\subset M$ separates $M$ 
(see Examples~\ref{expl:tight} and~\ref{expl:euler} 
in Section~\ref{sec:c-rd-surg}). 
 The case when the contact structure $\xi$ is overtwisted 
and the embedded torus $T\subset M$ does not separates $M$ 
is left to be discussed.

  Now, we discuss the case when we can not operate a contact round surgery 
directly. 
 Let $\xi$ be an overtwisted contact structure on a closed $3$-manifold $M$, 
and $T\subset (M,\xi)$ an embedded torus which does not separate $M$. 
 When the contact structure $\xi$ is overtwisted, 
there exists a case when $\releu T\ne 0$. 
 Assume that $\releu T\ne 0$. 

 First, we translate the condition in terms of the Euler characteristic 
of the regions divided by the dividing set. 
 Let $R_\pm\subset T$ be the positive and negative regions. 
 From Formulas~\eqref{eq:formulae}, we have $\releu T=2\chi(R_+)$ 
since $\chi(\Sigma)=\chi(T^2)=0$. 
 Therefore, the assumption $\releu T\ne0$ implies $\chi(R_+)\ne0$.  

\begin{figure}[htb]
  \centering
  \def\svgwidth{8cm}
  {\small 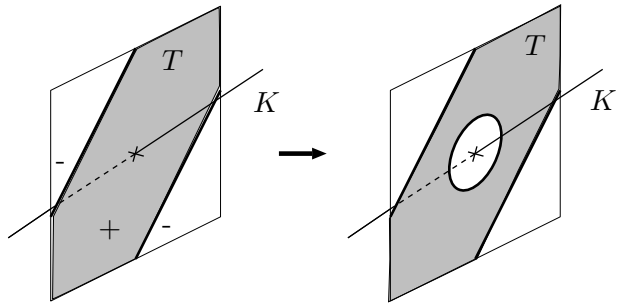}
  \caption{making a homotopically trivial dividing curve}
  \label{fig:mkdivcrv}
\end{figure}

  The Euler characteristic $\chi(R_+)$ of the positive region $R_+$ 
can be changed by modifying the contact structure $\xi$ as follows. 
 Recall that the embedded torus $T\subset(M,\xi)$ 
does not separate the manifold $M$. 
 Then we have a transverse knot $K\subset(M,\xi)$ 
which intersects the positive region $R_+\subset T$ once transversely. 
 By the Lutz twist along $K$ sufficiently close to $K$, 
we obtain a contact structure $\bar\xi$ modified from $\xi$ around $K$. 
 With respect to this $\bar\xi$, 
we have another homotopically trivial dividing curve on $T$ 
around the point where $K$ intersects $T$ (see Figure~\ref{fig:mkdivcrv}). 
%
%
 Let $\bar R_\pm$ denote the new positive and negative regions on $T$ 
with respect to $\bar\xi$. 
 Then we have 
\begin{equation*}
  \chi(\bar R_+)=\chi(R_+)-1,\qquad \chi(\bar R_-)=\chi(R_-)+1. 
\end{equation*}
 Similarly, by the simple Lutz twist along a transverse knot $K'$ 
which intersects $R_-$, 
the Euler characteristics are changed as 
\begin{equation*}
  \chi(\bar R'_+)=\chi(R_+)+1,\qquad \chi(\bar R'_-)=\chi(R_-)-1, 
\end{equation*}
where $\bar R'_\pm$ are the positive and negative regions 
after the simple Lutz twist along $K'$. 
 Then, by applying the simple Lutz twists suitably, 
we obtain a contact structure $\zeta$ on $M$ 
with respect to which $\releu[\zeta]{T}=0$ holds. 
 Note that we have not changed the manifold $M$ and the embedded torus $T$ 
but a contact structure on $M$. 

  Then we can conclude that any given round surgery of index~$2$ 
can be operated as a contact round surgery of index~$2$, 
after changing the contact structure if necessary. 
 In fact, even if the given $2$ dimensional torus $T$ 
in a contact $3$-manifold $(M,\xi)$ satisfies $\releu T\ne 0$, 
we have another contact structure $\zeta$ on $M$ with $\releu[\zeta]T=0$, 
as above. 
 Then we can operate the given round surgery of index~$2$ 
of $M$ along $T\subset M$ 
as a contact round surgery of index~$2$ of $(M,\zeta)$ along $T$. 

  This completes the proof of Theorem~\ref{thma}. \qed

%
%
\begin{rmk}
  Although the argument above is sufficient 
to prove the existence of a contact structure on a given manifold, 
we remark one thing for the discussion in the following subsection. 
 In the construction above, 
we use the Lutz twists other than contact round surgeries. 
 However, it had been proved that the simple Lutz twist is described 
as a pair of contact round surgeries (see Theorem~\ref{thm:Lz}). 
 Therefore, by making detours, 
we obtain a sequence of round surgeries of index~$1$ and index~$2$ 
from the given contact $3$-sphere to any closed orientable $3$-manifold 
which can be operated as direct contact round surgeries. 
 In other words, we can construct a contact structure 
on any closed orientable $3$-manifold from a contact $3$-sphere 
only by contact round surgeries 
without changing the intermediate contact structures. 
\end{rmk}

\subsection{Construction of any contact structure} \label{sec:prThmb}
  Theorem~\ref{thmb} is proved in this section. 
 We show that any contact structure on any closed oriented $3$-manifold 
is constructed from the standard contact $3$-sphere by contact round surgeries 
of index~$1$ and index~$2$. 
 It is proved in the following 2 steps. 
 In the first step, we show that any closed contact $3$-manifold 
is constructed from a $3$-dimensional sphere with some contact structure 
by contact round surgeries. 
 In the second step, we show that any contact structure 
on a $3$-dimensional sphere 
is constructed from the standard contact structure by contact round surgeries. 

\smallskip

\noindent \textbf{Step 1:} First, we verify that the procedure, 
topological round surgery, is reversible. 
 Recall that a round surgery implies attaching a round handle to a cobordism. 
 If an $n$-dimensional manifold $N$ is obtained from a manifold $M$ 
by a round surgery of index~$k$, 
the boundary of $W:=(M\times I)+R_k$ consists of $M$ and $N$. 
 From the Poincar\'e duality trick, 
the cobordism $W$ between $M$ and $N$ can be considered 
as $W=(N\times I)+R_{n-k-1}$ (see \cite{asimov}). 
 This implies that $M$ is obtained from $N$ by a round surgery of index~$n-k-1$.
 Therefore, Theorem~\ref{thm:r-surg} implies 
the $3$-dimensional sphere $S^3$ is obtained 
from any connected closed orientable $3$-manifold 
by round surgeries of index~$1$ and index~$2$. 

  Next, we show that contact round surgeries are reversible. 
 In the operation of a contact round surgery of index~$1$, 
we remove two contact solid tori and reglue a contact thickened torus. 
 Then the operation for the surgered manifold 
removing the glued contact thickened torus 
and regluing the removed contact solid tori 
recover the original contact manifold. 
 The second operation is a contact round surgery of index~$2$ 
from the surgered manifold to the original one. 
 In fact, from the definition in Subsection~\ref{sec:1crdsrg}, 
the solid tori are tubular neighborhoods of transverse knots, 
and the thickened torus is the invariant tubular neighborhood of a convex torus.
 The same property holds for contact round surgery of index~$2$. 
 A contact round surgery of index~$2$ is operated 
after isotoping the given torus with Condition~\eqref{eq:eulcond}. 
 The torus can be isotoped to a convex torus with no homotopically trivial 
dividing curve (see Subsection~\ref{sec:2crdsrg}). 
 Therefore, a contact round surgery of index~$2$ is also an operation 
concerning tubular neighborhoods of transverse knots 
and an invariant tubular neighborhood of a convex torus. 
 Then it is reversible. 

 Then, in order to show that the given closed contact $3$-manifold 
is obtained from a contact $3$-sphere by contact round surgeries 
of index~$1$ and index~$2$, 
it is sufficient to show that a contact $3$-sphere is obtained 
from the given closed contact $3$-manifold 
by contact round surgeries of index~$1$ and index~$2$. 
  Now, we construct a $3$-dimensional sphere with a contact structure 
from the given closed contact $3$-manifold. 
 Let $(M,\xi)$ be the given closed contact $3$-manifold. 
 From the argument above, there exists a sequence of round surgeries 
which makes $M$ the $3$-dimensional sphere. 
 By the same argument as in Subsection~\ref{sec:prThma}, 
this sequence with some detours can be operated 
as direct contact round surgeries. 
 Thus we obtain a $3$-dimensional sphere with some contact structure. 

\smallskip

\noindent \textbf{Step 2:} We discuss contact round surgeries 
between the standard contact structure and any other contact structure 
on the $3$-dimensional sphere $S^3$. 
 From Theorem~\ref{thm:classifS3}, in each homotopy class as plane fields 
other than the class of the standard contact structure, 
there exists a unique overtwisted contact structure up to isotopy. 
 In the class of the standard structure, 
there exists another overtwisted contact structure. 
 Recall that a simple Lutz twist changes the homotopy class 
of the contact structure under consideration 
(see Subsection~\ref{sec:dfnlztw}). 
 Then, it is known that all overtwisted contact structures on $S^3$ 
which are not homotopic to the standard one 
are constructed from the standard structure by the Lutz twists 
(see \cite{lutz}, \cite{eliash20years}). 
 The overtwisted structure homotopic to the standard contact structure 
is obtained from the standard structure by a full Lutz twist 
which is a consecutive two simple Lutz twists along the same transverse knot. 
 As we show in Subsection~\ref{sec:cstrlztw}, 
the simple Lutz twist is realized by a sequence of contact round surgeries 
(Theorem~\ref{thm:Lz}). 
 Therefore, all contact structures on the $3$-dimensional sphere $S^3$ 
are obtained from the standard contact structure by contact round surgeries. 

  Thus Theorem~\ref{thmb} has been proved. \qed

%
%


\begin{thebibliography}{99999}
\bibitem[Ad1]{art19} J.~Adachi, 
  \textit{Contact round surgery and symplectic round handlebodies}, 
  Internat.\ J.\ Math.\ \textbf{25} (2014), 1450050, 25 pp. 

\bibitem[Ad2]{art20} J.~Adachi, 
  \textit{Contact round surgery and Lutz twists}, 
  (preprint). 

\bibitem[As]{asimov} D.~Asimov, 
  \textit{Round handles and non-singular Morse-Smale flows}, 
  Ann.\ of Math.\ (2) \textbf{102} (1975), 41--54. 

\bibitem[Ba]{baykur} R.~\.I.~Baykur, 
  \textit{Topology of broken Lefschetz fibrations 
    and near-symplectic four-manifolds}, 
  Pacific J.\ Math. \textbf{240} (2009), 201--230. 

\bibitem[BoElMu]{bem} M.~Borman, Ya.~Eliashberg, E.~Murphy, 
  \textit{Existence and classification of overtwisted contact structures 
    in all dimensions}, 
  Acta Math. \textbf{215} (2015), 281--361. 

\bibitem[Br]{bredon} G.~E.~Bredon, 
  \textit{Topology and geometry}, 
  Graduate Texts in Mathematics, \textbf{139}, 
  Springer-Verlag, New York, 1993. 

\bibitem[DGe]{dingei04} F.~Ding and H.~Geiges, 
  \textit{A Legendrian surgery presentation of contact 3-manifolds}, 
  Math.\ Proc.\ Cambridge Philos.\ Soc.\ \textbf{136} (2004), 583--598. 

\bibitem[El1]{eliashOT} Ya.~Eliashberg, 
  \textit{Classification of overtwisted contact structures on $3$-manifolds}, 
  Invent.\ Math.\ \textbf{98} (1989), 623--637. 

\bibitem[El2]{eliashCharStein} Ya.~Eliashberg, 
  \textit{Topological characterization of Stein manifolds of dimension $>2$}, 
  Internat.\ J.\ Math. \textbf{1} (1990), 29--46. 

\bibitem[El3]{eliash20years} Ya.~Eliashberg, 
  \textit{Contact $3$-manifolds twenty years since J.~Martinet's work}, 
  Ann.\ Inst.\ Fourier (Grenoble) \textbf{42} (1992), 165--192. 

\bibitem[Et]{etsurg} J.~Etnyre, 
  \textit{On contact surgery}, 
  Proc.\ Amer.\ Math.\ Soc.\ \textbf{136} (2008), 3355--3362. 

\bibitem[EtGh]{etgh} J.~Etnyre and R.~Ghrist, 
  \textit{Gradient flows within plane fields}, 
  Comment.\ Math.\ Helv.\ \textbf{74} (1999), 507--529.

\bibitem[Ga]{gay06}  Gay, David T. 
  \textit{Four-dimensional symplectic cobordisms containing three-handles}, 
  Geom.\ Topol.\ \textbf{10} (2006), 1749--1759. 

\bibitem[Ge]{geitext} H.~Geiges, 
  An introduction to contact topology, 
  \emph{Cambridge Studies in Advanced Mathematics} \textbf{109}, 
  Cambridge University Press, Cambridge, 2008. 

\bibitem[Gi1]{girouxConv} E.~Giroux, 
  \textit{Convexit\'e en topologie de contact}, 
  Comment.\ Math.\ Helv. \textbf{66} (1991), 637--677. 

\bibitem[Gi2]{girouxTor} E.~Giroux, 
  \textit{Une structure de contact, m\^eme tendue, est plus ou moins tordue}, 
  Ann.\ Sci.\ \'Ecole Norm.\ Sup.\ (4)  \textbf{27} (1994), 697--705. 

\bibitem[Gi3]{giroux99} E.~Giroux, 
  \textit{Une infinit\'e de structures de contact tendues sur une infinit\'e 
    de vari\'et\'es}, 
  Invent.\ Math.\ \textbf{135} (1999), 789--802. 

\bibitem[Gi4]{girouxCrit} E.~Giroux, 
  \textit{Structures de contact sur les vari\'et\'es fibr\'ees 
    en cercles audessus d'une surface},   
  Comment.\ Math.\ Helv.\ \textbf{76} (2001), 218--262. 

\bibitem[Ho]{hondaI} K.~Honda, 
  \textit{On the classification of tight contact structures. I}, 
  Geom.\ Topol.\ \textbf{4} (2000), 309--368. 

\bibitem[Hu]{huang} Y.~Huang, 
  \textit{A proof of the classification theorem 
    of overtwisted contact structures via convex surface theory}, 
  J.\ Symplectic Geom.\ \textbf{11} (2013), 563--601. 

\bibitem[K]{kanda} Y.~Kanda, 
  \textit{The classification of tight contact structures on the $3$-torus}, 
  Comm.\ Anal.\ Geom. \textbf{5} (1997), 413--438. 

\bibitem[L]{lutz} R.~Lutz, 
  \textit{Structures de contact sur les fibr\'es principaux en cercles 
    de dimension trois}, 
  Ann.\ Inst.\ Fourier (Grenoble) \textbf{27} (1977), 1--15. 

\bibitem[Mar]{martinet} J.~Martinet, 
  \textit{Formes de contact sur les vari\'et\'es de dimension $3$}, 
  Proceedings of Liverpool Singularities Symposium, II (1969/1970),  
  pp.~142--163, Lecture Notes in Math., Vol.~\textbf{209}, 
  Springer, Berlin, 1971. 

\bibitem[MasNWen]{mnw} Massot, Patrick; Niederkr\"uger, Klaus; Wendl, Chris 
  \textit{Weak and strong fillability of higher dimensional contact manifolds}, 
  Invent.\ Math. \textbf{192} (2013), 287--373. 

\bibitem[Mi]{miyoshi} S.~Miyoshi, 
  \textit{Foliated round surgery of codimension-one foliated manifolds}, 
  Topology \textbf{21} (1982), 245--261. 

\bibitem[Mo]{morgan} J.~Morgan, 
  \textit{Nonsingular Morse-Smale flows on $3$-dimensional manifolds}, 
  Topology \textbf{18} (1979), 41--53. 

\bibitem[OzbSt]{osbook} B.~Ozbagci and A.~Stipsicz, 
  Surgery on contact 3-manifolds and Stein surfaces, 
  \textit{Bolyai Society Mathematical Studies}, \textbf{13}, 
  Springer-Verlag, Berlin, 2004. 

\bibitem[OzsSz]{ozsz} P.~Ozsv\'ath and Z.~Szab\'o, 
  \textit{Heegaard Floer homology and contact structures}, 
  Duke Math.\ J.\ \textbf{129} (2005), 39--61.

\bibitem[St]{stipsicz} A.~I.~Stipsicz, 
  \textit{Surgery diagrams and open book decompositions of contact 3-manifolds},
  Acta Math.\ Hungar.\  \textbf{108} (2005),  71--86.


\bibitem[V]{vogel} T.~Vogel, 
  \textit{Existence of Engel structures}, 
  Ann.\ of Math.\ (2) \textbf{169} (2009), 79--137. 

\bibitem[Wei]{weinsteinHdl} A.~Weinstein, 
  \textit{Contact surgery and symplectic handlebodies}, 
  Hokkaido Math.\ J.\ \textbf{20} (1991), 241--251. 
\end{thebibliography}
\end{document}